\documentclass{revtex4}
\usepackage{graphicx}
\usepackage{epstopdf}
\usepackage{pstricks}
\usepackage{amsmath,amssymb}
\usepackage{xcolor}

\begin{document}
\title{A derivation of boundary conditions for the space-fractional operator with order $0<\alpha<1$  at the vertices of metric graphs}
\author{Sabirov K.K.$^{1,2}$, Ashurov R.R.$^1$}
\affiliation{$^1$Institute of Mathematics named after V.I.Romanovsky of the Academy of Sciences of the Republic of Uzbekistan, University street 9, Tashkent, 100174, Uzbekistan\\
$^2$Tashkent University of  Information Technology, Amir Temur
street 108, Tashkent 100200, Uzbekistan}
\begin{abstract}
We consider the space-fractional operator with order $0<\alpha<1$ on the metric star graph. The boundary conditions at the vertices of the metric star graph providing the self-adjointness of the operator are derived. The obtained result is extended to the other topologies of the metric graphs.
\end{abstract}
\maketitle
\section{Introduction}

Fractional differential operators have attracted  much  attention in the last few decades due to their applications in numerous seemingly distant ﬁelds of science and technology (see \cite{Hilfer}-\cite{Umarov}). The spectral analysis of the nonself-adjoint fractional differential operators was studied in detail by T.S.Aleroev \cite{Aleroev1999}-\cite{Aleroev2016} and the self-adjoint boundary conditions at the vertices of the metric graphs for the operators degenerated by Schr\"odinger, Dirac and Bogoliubov de Gennes equations are studied in \cite{Schrader}-\cite{Sabirov}.

In this paper we consider the fractional differential operator on metric graphs. This paper is organized as follows. There are given the mathematical preliminaries as  fractional integrals and derivatives in the section 2. In the section 3 we consider the fractional differential operator with order $0<\alpha<1$ on the star graph and we derive the self-adjoint boundary conditions at vertices of the star graph. In the section 4 we have given the extensions to the other topologies of the obtained result and finally, the last section presents conclusions.

\section{Mathematical preliminaries}

{\bf Definition 1:} Let $\Omega=[a;b]$ be finite interval on real axis $R$. The Riemann-Lioville fractional integrals denoted by $I_{a+}^\alpha$ and $I_{b-}^\alpha$ of order $\alpha\in C\, (Re(\alpha)>0)$ are defined as
\begin{equation}
\left(I_{a+}^\alpha f\right)(x):=\frac{1}{\Gamma(\alpha)}\underset{a}{\overset{x}{\int}}\frac{f(t)dt}{(x-t)^{1-\alpha}},\,(x>a;\,Re(\alpha>0))\label{rli1}
\end{equation}
and
\begin{equation}
\left(I_{b-}^\alpha f\right)(x):=\frac{1}{\Gamma(\alpha)}\underset{x}{\overset{b}{\int}}\frac{f(t)dt}{(t-x)^{1-\alpha}},\,(x<b;\,Re(\alpha>0))\label{rli2}
\end{equation}
respectively, where $\Gamma(\alpha)$ is the Gamma function, integrals (\ref{rli1}) and (\ref{rli2}) are called the left-sided and the right-sided fractional integrals.

{\bf Definition 2} The Riemann-Lioville derivatives denoted $D_{a+}^\alpha y$ and $D_{b-}^\alpha y$ of order $\alpha\in {\bf C}$ $(Re(\alpha)\geq 0)$ are defined by
\begin{equation}
\left(D_{a+}^\alpha y\right)(x):=\frac{1}{\Gamma(1-\alpha)}\left(\frac{d}{dx}\right)^n\underset{a}{\overset{x}{\int}}\frac{y(t)dt}{(x-t)^{\alpha-n+1}},\,(n=\left[Re(\alpha)\right]+1,\,x>a),\label{rlfd1}
\end{equation}
and
\begin{equation}
\left(D_{b-}^\alpha y\right)(x):=\frac{1}{\Gamma(1-\alpha)}\left(-\frac{d}{dx}\right)^n\underset{x}{\overset{b}{\int}}\frac{y(t)dt}{(t-x)^{\alpha-n+1}},\,(n=\left[Re(\alpha)\right]+1,\,x<b),\label{rlfd2}
\end{equation}
respectively, where $[Re(\alpha)]$ means integer part of $Re(\alpha)$. In particular, $0<\alpha<1$, then
\begin{equation}
\left(D_{a+}^\alpha y\right)(x)=\frac{1}{\Gamma(1-\alpha)}\frac{d}{dx}\underset{a}{\overset{x}{\int}}\frac{y(t)dt}{(x-t)^\alpha},\,(x>a),\label{rlfdf1}
\end{equation}
and
\begin{equation}
\left(D_{b-}^\alpha y\right)(x)=-\frac{1}{\Gamma(1-\alpha)}\frac{d}{dx}\underset{x}{\overset{b}{\int}}\frac{y(t)dt}{(t-x)^\alpha},\,(x<b)\label{rlfdf2}
\end{equation}

We introduce the following corollary.

{\bf Corollary 1} (see \cite{Kilbas}):
If $0<\alpha<1$ and $y(x)\in AC[a,b]$, then
\begin{equation}
    \left(D_{a+}^\alpha y\right)(x)=\frac{1}{\Gamma(1-\alpha)}\left[\frac{y(a)}{(x-a)^\alpha}+\underset{a}{\overset{x}{\int}}\frac{y'(t)dt}{(x-t)^\alpha}\right],\label{rlcfd1}
\end{equation}
and
\begin{equation}
    \left(D_{b-}^\alpha y\right)(x)=\frac{1}{\Gamma(1-\alpha)}\left[\frac{y(b)}{(b-x)^\alpha}-\underset{x}{\overset{b}{\int}}\frac{y'(t)dt}{(t-x)^\alpha}\right].\label{rlcfd2}
\end{equation}
\section{Formulation of the problem}

\begin{figure}[ht!]
\centering
\includegraphics[scale=0.4]{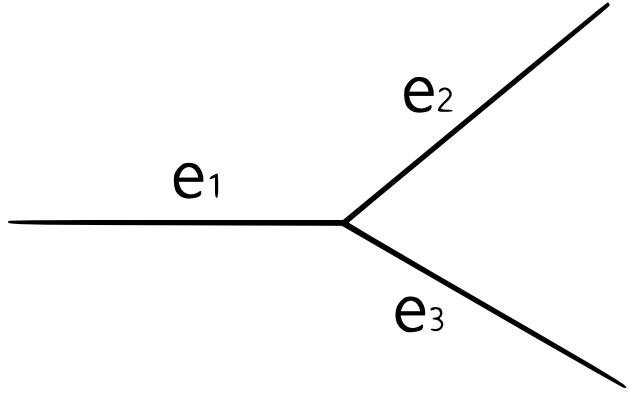}
\caption{The metric star graph} \label{pic1}
\end{figure}

We consider the star graph $G$ with tree bonds $e_j,\,j=1,2,3$, for which a coordinate $x_j$ is assigned. Choosing the origin of coordinates at the vertex, 0, for bond $e_j$ we put $x_j\in[0;l_j]]$. Let the Hilbert space be given as
\begin{eqnarray}
{\bf L_2}(G):=\underset{j=1}{\overset{3}{\oplus}}L_2([0;l_j])\label{l2},
\end{eqnarray}
and the space of the absolutely continuous functions be given as
\begin{equation}
{\bf AC}(G):=\underset{j=1}{\overset{3}{\oplus}}AC([0;l_j]).\label{ac}
\end{equation}

We use the shorthand notation $\psi_{j}(x)$ for $\psi_j(x_j)$ where $x$ is the coordinate on the bond $j$ to which the component $\psi_j$ refers.
The scalar product in ${\bf L_2}(G)$ is the following
\begin{equation}
\langle\varphi,\psi\rangle:=\underset{j=1}{\overset{3}{\sum}}\underset{e_j}{\int}\varphi_j(x)\cdot\psi_j^*(x)dx.\label{sp1}
\end{equation}
We consider the following fractional operator with order $\alpha$
\begin{equation}
A_0:=(A_0^{(1)},A_0^{(2)},A_0^{(3)}),\label{fop1}
\end{equation}
with the domain of the definition $D_0:=\left\{\varphi|\,\varphi=(\varphi_1,\varphi_2,\varphi_3)\in{\bf AC}(G),\,\varphi_j(0)=\varphi_j(l_j)=0,\,j=1,2,3\right\}\subseteq{\bf AC}(G)\subseteq{\bf L_2}(G)$, where $A_0^{(j)}:=D_{0+}^\alpha+D_{l_j-}^\alpha\,(0<\alpha<1)$.

{\bf Lemma 1:} The operator $A_0$ is symmetric operator.

Proof. We take $\varphi\in D_0$ and
\begin{eqnarray}
\langle A_0\varphi,\psi\rangle=\underset{j=1}{\overset{3}{\sum}}\langle A_0^{(j)}\varphi_j,\psi_j\rangle=\underset{j=1}{\overset{3}{\sum}}\int_{0}^{l_j}\left(D_{0+}^\alpha+D_{l_j-}^\alpha\right)\varphi_j(x)\cdot\psi_j^*(x)dx=\nonumber\\
=\underset{j=1}{\overset{3}{\sum}}\int_{0}^{l_j}\frac{1}{\Gamma(1-\alpha)}\left(\frac{d}{dx}\int_{0}^{x}\frac{\varphi_j(t)}{(x-t)^\alpha}dt-\frac{d}{dx}\int_{x}^{l_j}\frac{\varphi_j(t)}{(t-x)^\alpha}dt\right)\psi_j^*=\nonumber\\
=\underset{j=1}{\overset{3}{\sum}}\left[\psi_j^*(l_j)\cdot\left(I_{0+}^{1-\alpha}\varphi_j\right)(l_j-0)-\psi_j^*(0)\cdot\left(I_{0+}^{1-\alpha}\varphi_j\right)(+0)\right.-\nonumber\\
-\left.\psi_j^*(l_j)\cdot\left(I_{l_j-}^{1-\alpha}\varphi_j\right)(l_j-0)+\psi_j^*(0)\cdot\left(I_{l_j-}^{1-\alpha}\varphi_j\right)(+0)\right]-\nonumber\\
-\underset{j=1}{\overset{3}{\sum}}\frac{1}{\Gamma(1-\alpha)}\int_{0}^{l_j}\varphi_j(t)\cdot\left[\int_{t}^{l_j}\frac{\psi_j^{*'}(x)}{(x-t)^\alpha}dx-\int_{0}^{t}\frac{\psi_j^{*'}(x)}{(t-x)^\alpha}dx\right]dt=\nonumber\\
=-\underset{j=1}{\overset{3}{\sum}}\left[\psi_j^*(0)\cdot\left(I_{0+}^{1-\alpha}\varphi_j\right)(+0)+\psi_j^*(l_j)\cdot\left(I_{l_j-}^{1-\alpha}\varphi_j\right)(l_j-0)\right]+\nonumber\\
+\underset{j=1}{\overset{3}{\sum}}\int_{0}^{l_j}\varphi_j(t)\cdot\left(\left(D_{0+}^\alpha+D_{l_j-}^\alpha\right)\psi_j(t)\right)^*dt=\langle\varphi,A_0^*\psi\rangle,
\end{eqnarray}
where $A_0\subset A_0^*$.
Lemma 1 is proved.

For the self-adjoint extension of operator $A_0$ (with the domain of the definition $D_0$) we consider the following bi-linear form (skew-Hermitian)
\begin{equation}
\Omega(\varphi,\psi)=\langle A\varphi,\psi\rangle-\langle\varphi,A\psi\rangle\label{skewh1}
\end{equation}
with a differential expression $A:=(D_{0+}^\alpha+D_{l_1-}^\alpha,D_{0+}^\alpha+D_{l_2-}^\alpha,D_{0+}^\alpha+D_{l_3-}^\alpha)$ on ${\bf AC}(G)\subseteq{\bf L_2}(G)$ $\,(0<\alpha<1)$

From $\Omega(\varphi,\psi)=0$ we can obtain boundary conditions at the vertices of the graph which provide the self-adjoint extension of the operator $A_0$.

{\bf Theorem 1:} The following boundary conditions provide the self-adjoint extension of the operator $A_0$ on ${\bf AC}(G)\subseteq{L_2}(G)$
\begin{eqnarray}
\alpha_1\varphi_1|_{x=0}=\alpha_2\varphi|_{x=0}=\alpha_3\varphi|_{x=0},\nonumber\\
\frac{1}{\alpha_1}\left(I_{+0}^{1-\alpha}\varphi_1\right)(+0)+\frac{1}{\alpha_2}\left(I_{+0}^{1-\alpha}\varphi_2\right)(+0)+\frac{1}{\alpha_3}\left(I_{+0}^{1-\alpha}\varphi_3\right)(+0)=0,\nonumber\\
\varphi_j(l_j)=0,\,j=1,2,3,\label{bc1}
\end{eqnarray}
where $\alpha_j,\,j=1,2,3$ are real nonzero constants.

Proof of theorem 1. 
\begin{eqnarray}
&\Omega(\varphi,\psi)=\langle A\varphi,\psi\rangle-\langle\varphi,A\psi\rangle=\nonumber\\
&=-\underset{j=1}{\overset{3}{\sum}}\left[\psi_j^*(l_j)\cdot\left(I_{l_j-}^{1-\alpha}\varphi_j\right)(l_j-0)+\psi_j^*(0)\cdot\left(I_{0+}^{1-\alpha}\psi_j^*\right)(+0)\right]=0.\label{skewh2}
\end{eqnarray}
The theorem is proved.

{\bf Theorem 2:} The ordinary fractional differential equation
\begin{eqnarray}
A_j\varphi_j=k_j\varphi_j\label{eq1}
\end{eqnarray}
with boundary conditions (\ref{bc1}) has a solution.

{\bf Proof.} One can derive the following Volterra integral equation from Eq. (\ref{eq1})
\begin{equation}
\varphi_j(x)=\frac{b_j}{\Gamma(\alpha)}x^{\alpha-1}-\frac{1}{\Gamma(\alpha)}\underset{0}{\overset{x}{\int}}\frac{-k_j\varphi_j(s)+\left(D_{l_j-}^{\alpha}\varphi_j\right)(s)}{\left(x-s\right)^{1-\alpha}}ds,\label{Volt1}
\end{equation}
with
\begin{equation}
\frac{b_1}{\alpha_1}+\frac{b_2}{\alpha_2}+\frac{b_3}{\alpha_3}=0.\label{constrain1}
\end{equation}
The last equation we put into $\varphi_j(l_j)=0$ and find the following
\begin{eqnarray}
\frac{b_j\cdot l_j^{\alpha-1}}{\Gamma(\alpha)}-\frac{1}{\Gamma(\alpha)}\underset{0}{\overset{l_j}{\int}}\frac{-k_j\varphi_j(s)+\left(D_{l_j-}^{\alpha}\varphi_j\right)(s)}{\left(l_j-s\right)^{1-\alpha}}ds=0,\nonumber
\end{eqnarray}
or
\begin{eqnarray}
\underset{0}{\overset{l_j}{\int}}\frac{b_j\cdot l_j^{\alpha-1}(l_j-s)^{1-\alpha}+k_j\varphi_j(s)-\left(D_{l_j-}^{\alpha}\varphi_j\right)(s)}{\left(l_j-s\right)^{1-\alpha}}ds=0.\nonumber
\end{eqnarray}
From the last equation we can get
\begin{eqnarray}
\left(D_{l_j-}^{\alpha}\varphi_j\right)(s)-k_j\varphi_j(s)=b_j\cdot l_j^{\alpha-2}(l_j-s)^{1-\alpha}.\label{eq2}
\end{eqnarray}
A solution of this equation is the following
\begin{eqnarray}
\varphi_j(s)=c_j\cdot(l_j-s)^{\alpha-1}E_{\alpha,\alpha}\left[k_j(l_j-s)^\alpha\right]+b_jl_j^{\alpha-2}\underset{s}{\overset{l_j}{\int}}(t-s)^{\alpha-1}E_{\alpha,\alpha}\left[k_j(t-s)^\alpha\right](l_j-t)^{1-\alpha}dt,\label{sol1}
\end{eqnarray}
where $c_j$ are real constants and $E_{\alpha,\beta}(x)$ is Mittag-Leffler function. Using boundary conditions (\ref{bc1}) at the vertices we have the following constrains
\begin{eqnarray}
\frac{\alpha_{2,3}}{\alpha_1}=\frac{b_1}{b_{2,3}}\cdot\left(\frac{l_1}{l_{2,3}}\right)^{\alpha-1},\,\alpha_1\cdot c_1l_1^{\alpha-1}=\alpha_2\cdot c_2l_2^{\alpha-1}=\alpha_3\cdot c_3l_3^{\alpha-1},\,k_1l_1^\alpha=k_2l_2^\alpha=k_3l_3^\alpha.\label{constrain2}
\end{eqnarray}
From this constrains (\ref{constrain1}) and (\ref{constrain2}) we can get
\begin{eqnarray}
\frac{b_1}{l_1k_1}+\frac{b_2}{l_2k_2}+\frac{b_3}{l_3k_3}=0.\label{constrain3}
\end{eqnarray}
The theorem 2 is proved.

Moreover, the vertex boundary conditions for the self-adjoint extension of the operator $A_0:=(A_0^{(1)},A_0^{(2)},...,A_0^{(N)})$ (with the domain of definition $D_0:=\left\{\varphi|\,\varphi=(\varphi_1,\varphi_2,...,\varphi_N)\in{\bf AC}(G),\,\varphi_j(0)=\varphi_j(l_j)=0,\,j=1,2,...,N\right\}\subseteq{\bf AC}(G)\subseteq{\bf L_2}(G)$, where $A_0^{(j)}:=D_{0+}^\alpha+D_{l_j-}^\alpha$ $(0<\alpha<1)$) on the star graph $G$ with $N$ bonds $e_j,\,j=1,2,...,N$ are the following
\begin{eqnarray}
\alpha_1\psi_1|_{x=0}=\alpha_2\psi_2|_{x=0}=...=\alpha_N\psi_N|_{x=0},\nonumber\\
\frac{1}{\alpha_1}\left(I_{+0}^{1-\alpha}\psi_1\right)(+0)+\frac{1}{\alpha_2}\left(I_{+0}^{1-\alpha}\psi_2\right)(+0)+...+\frac{1}{\alpha_N}\left(I_{+0}^{1-\alpha}\psi_N\right)(+0)=0,\nonumber\\
\psi_j(l_j)=0,\,j=1,2,...,N.\label{bc2}
\end{eqnarray}

\section{Extensions to other topologies}

\begin{figure}[ht!]
\centering
\includegraphics[scale=0.3]{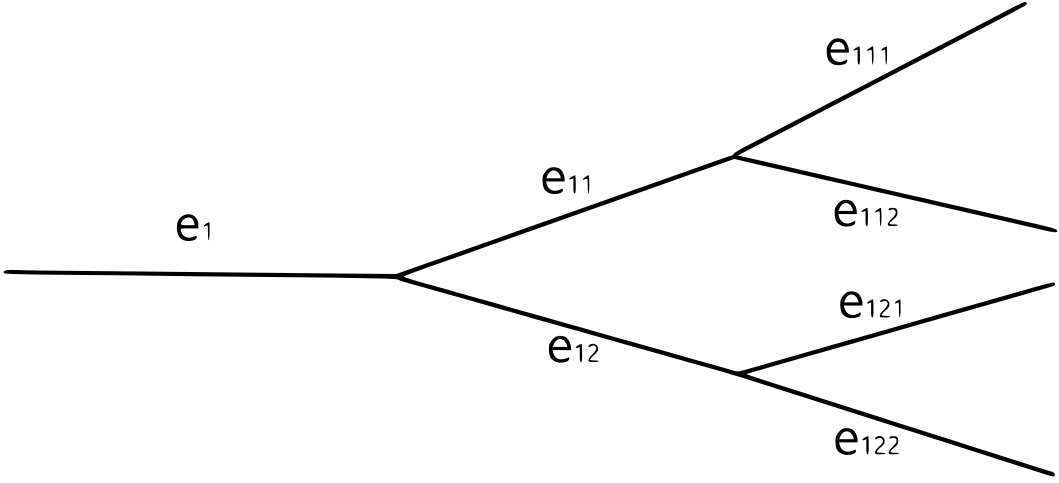}
\caption{The metric tree graph} \label{pic2}
\end{figure}
We consider a tree graph with bonds $e_1\sim[0;l_1],\,e_{1i}\sim[0;l_{1i}],\,e_{1ij}\sim[0;l_{1ij}],\,i,j=1,2$ plotted in Fig.\ref{pic2}.

{\bf Theorem 3.} The boundary conditions at the vertices for the self-adjoint extension of the operator $A_0:=(A_0^{(1)},A_0^{(11)},A_0^{(12)},A_0^{(111)},A_0^{(112)},A_0^{(121)},A_0^{(122)})$ (with the domain of definition $D_0:=\left\{\varphi|\,\varphi=(\varphi_1,\varphi_{11},\varphi_{12},\varphi_{111},\varphi_{112},\varphi_{121},\varphi_{122})\in{\bf AC}(G),\,\varphi_1(0)=\varphi_1(l_1)=\varphi_{1i}(0)=\varphi_{1i}(l_{1i})=\varphi_{1ij}(0)=\varphi_{1ij}(l_{1ij})=0,\right.$ $\left.i,j=1,2\right\}\subseteq{AC}(G)\subseteq{\bf L_2}(G)$, where $A_0^{(1)}:=D_{0+}^\alpha+D_{l_1-}^\alpha,\,A_0^{(1i)}:=D_{0+}^\alpha+D_{l_{1i}-}^\alpha,\,A_0^{(1ij)}:=D_{0+}^\alpha+D_{l_{1ij}-}^\alpha,\,i,j=1,2$ $(0<\alpha<1)$) on the tree graph $G$ with bonds $e_1,\,e_{1i},\,e_{1ij},\,i,j=1,2$ are the following
\begin{eqnarray}
\alpha_1\psi_1|_{x=l_1}=\alpha_{11}\psi_{11}|_{x=0}=\alpha_{12}\psi_{12}|_{x=0},\nonumber\\
\frac{1}{\alpha_1}\left(I_{l_1-}^{1-\alpha}\psi_1\right)(l_1-)+\frac{1}{\alpha_{11}}\left(I_{+0}^{1-\alpha}\psi_{11}\right)(+0)+\frac{1}{\alpha_{12}}\left(I_{+0}^{1-\alpha}\psi_{12}\right)(+0)=0,\nonumber\\
\alpha'_{1i}\psi_{1i}|_{x=l_{1i}}=\alpha_{1i1}\psi_{1i1}|_{x=0}=\alpha_{1i2}\psi_{1i2}|_{x=0},\,i=1,2,\nonumber\\
\frac{1}{\alpha'_{1i}}\left(I_{l_{1i}-}^{1-\alpha}\psi_{1i}\right)(l_{1i}-)+\frac{1}{\alpha_{1i1}}\left(I_{+0}^{1-\alpha}\psi_{1i1}\right)(+0)+\frac{1}{\alpha_{1i2}}\left(I_{+0}^{1-\alpha}\psi_{1i2}\right)(+0)=0,\,i=1,2,\nonumber\\
\psi_1(0)=\psi_{1ij}(l_{1ij})=0,\,i,j=1,2.\label{bc3}
\end{eqnarray}

\begin{figure}[ht!]
\centering
\includegraphics[scale=0.3]{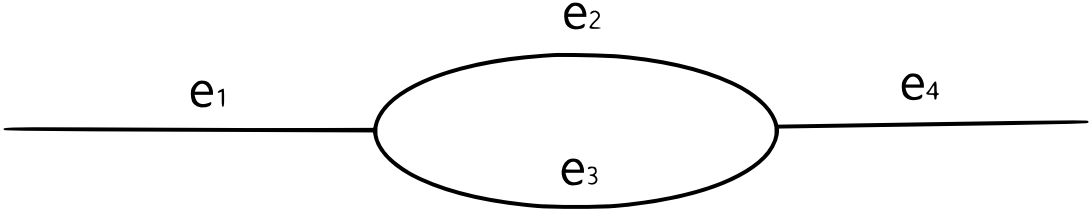}
\caption{The metric loop graph} \label{pic3}
\end{figure}

Now we consider a loop graph with bonds $e_j\sim[0;l_j],\,j=1,2,3,4$ plotted in Fig.\ref{pic3}.

{\bf Theorem 4.} The boundary conditions at the vertices for the self-adjoint extension of the operator $A_0:=(A_0^{(1)},A_0^{(2)},A_0^{(3)},A_0^{(4)})$ (with the domain of definition $D_0:=\left\{\varphi|\,\varphi=(\varphi_1,\varphi_{2},\varphi_{3},\varphi_{4})\in{\bf AC}(G),\,\varphi_j(0)=\varphi_j(l_j)=0,\,j=1,2,3,4\right\}\subseteq{AC}(G)\subseteq{\bf L_2}(G)$, where $A_0^{(j)}:=D_{0+}^\alpha+D_{l_j-}^\alpha,\,j=1,2,3,4$ $(0<\alpha<1)$) on the tree graph $G$ with bonds $e_j,\,j=1,2,3,4$ are the following
\begin{eqnarray}
\alpha_1\psi_1|_{x=l_1}=\alpha_{2}\psi_{2}|_{x=0}=\alpha_{3}\psi_{3}|_{x=0},\nonumber\\
\frac{1}{\alpha_1}\left(I_{l_1-}^{1-\alpha}\psi_1\right)(l_1-)+\frac{1}{\alpha_{2}}\left(I_{+0}^{1-\alpha}\psi_{2}\right)(+0)+\frac{1}{\alpha_{3}}\left(I_{+0}^{1-\alpha}\psi_{3}\right)(+0)=0,\nonumber\\
\alpha'_{2}\psi_{2}|_{x=l_{2}}=\alpha'_{3}\psi_{3}|_{x=l_3}=\alpha_{4}\psi_{4}|_{x=0},\,i=1,2,\nonumber\\
\frac{1}{\alpha'_{2}}\left(I_{l_{2}-}^{1-\alpha}\psi_{2}\right)(l_{2}-)+\frac{1}{\alpha'_{3}}\left(I_{l_3-}^{1-\alpha}\psi_{3}\right)(l_3-)+\frac{1}{\alpha_{4}}\left(I_{+0}^{1-\alpha}\psi_{4}\right)(+0)=0,\,i=1,2,\nonumber\\
\psi_1(0)=\psi_{4}(l_{4})=0.\label{bc4}
\end{eqnarray}

Proof of the theorems 3 and 4 is similarly with proof of the theorem 1.
\section{Conclusions}

 In this paper we have studied the spectral analysis of the fractional differential operator with order $0<\alpha<1$ on the metric star graph. We have obtained the self-adjoint boundary conditions at the branched point. The obtained result has been extended to the other topologies as metric tree and loop graphs.


\begin{thebibliography}{99}
\bibitem{Hilfer} R. Hilfer, Applications of fractional calculus in physics, Stuttgart, Germany, 2000.
\bibitem{Kilbas} Anatoly A. Kilbas, Hari M. Srivastava, Juan J. Trujillo, Theory and applications of fractional differential equations, Amsterdam, the Netherlands, 2006.
\bibitem{Umarov} S. Umarov, Introduction to Fractional and Pseudo-Diﬀerential Equations with Singular Symbols, Springer, 2015.
\bibitem{Schrader} V. Kostrykin, R. Schrader, Kirchhoff’s rule for quantum wires, Journal of Physics A: Mathematical and General, {\bf 32} (1999), 595–630.
\bibitem{Bolte} J. Bolte, J. Harrison, Spectral statistics for the Dirac operator on graphs, Journal of Physics A: Mathematical and General, {\bf 36} (2003), 2747–2769.
\bibitem{Harrison} J. M. Harrison, T. Weyand, K. Kirsten, Zeta functions of the Dirac operator on quantum graphs, Journal of Mathematical Physics, {\bf 57} (2016), 102301.
\bibitem{Sabirov} K. Sabirov, D. Jumanazarov, J. Yusupov, D. Matrasulov, Bogoliubov de Gennes equation on metric graphs, Physics Letters A, {\bf 382} (2018), 2856–2860.
\bibitem{Aleroev1999} T.  S.  Aleroev, A  Problem  on  the  Zeros  of the  Mittag-Lefiter  Function and  the  Spectrum  of a  Fractional-Order Differential  Operator, Differential  Equations,  {\bf 36},  9 (2000),  pp.  1414-1415.
\bibitem{Aleroev2014} T. S. Aleroev, Kh. T. Aleroeva, On a Class of Fractional Non-Selfadjoint Operators Associated with Diﬀerential Equations, Russian Mathematics (Iz. VUZ),  {\bf 58}, 10 (2014), pp. 1–9.
\bibitem{Aleroev2016} T. S. Aleroev, H. T. Aleroeva, On one class of nonself-adjoint operators associated with diﬀerential equations of fractional order, Journal of Mathematical Sciences, {\bf 214}, 2 (2016).
\end{thebibliography}
\end{document}